\documentclass[12pt]{article}   

\usepackage{a4,psfig,amssymb,verbatim,amsmath}

\begin{document}
\newcommand{\ddt}{\partial \over \partial t}
\newcommand{\ddx}{\partial \over \partial x}
\newcommand{\ddu}{\partial \over \partial u}
\def\dhp#1{ \mathop{#1}\limits_{h}}
\def\dphh#1{ \mathop{#1}\limits_{h \bar h}}
\def\da#1{ \mathop{#1}\limits_{+\tau}}
\def\db#1{ \mathop{#1}\limits_{-\tau}}
\def\dc#1{ \mathop{#1}\limits_{\pm \tau}}
\def\dd#1{ \mathop{#1}\limits_{+h}}
\def\df#1{ \mathop{#1}\limits_{-h}}
\def\dpm#1{ \mathop{#1}\limits_{\pm h}}
\def\dg#1{ \mathop{#1}\limits_{\pm h}}
\def\dh#1{ \mathop{#1}\limits_ h}

\begin{center}
{\Large {\bf Symmetry--preserving discrete schemes for \\
\medskip
some heat transfer equations}.}
\end{center}

\bigskip

\begin{center}

{\large 
 {Bakirova, M.I.}$^{*}$,  Dorodnitsyn, V.A.$^{*}$ and 
Kozlov, R.V.$^{\dag}$ } 

\bigskip
\bigskip

${}^{*}$
Keldysh Institute of Applied Mathematics of Russian Academy of Science, \\
Miusskaya Pl.~4, Moscow, 125047, Russia\\

${}^{\dag}$  Department of Mathematical Sciences, NTNU, N-7491,
Trondheim, Norway \\

\end{center}

\vspace{1.0 ex}


\bigskip

\begin {quotation}
\noindent {\bf Abstract.} Lie group analysis of differential 
equations is a
generally recognized method, which provides invariant solutions,
integrability, conservation laws etc. 
In this paper we present three
characteristic  examples of the construction of invariant difference
equations and meshes, where the original continuous symmetries are
preserved in discrete models. 
Conservation of symmetries 
in difference modeling helps to retain qualitative properties of 
the differential 
equations in their difference counterparts.

\end{quotation}

\bigskip

\noindent {\bf 1. Introduction}

\medskip

Symmetries are intrinsic and fundamental features of the differential 
equations of mathematical physics. Consequently, they should be 
retained when discrete analogs of such equations are constructed. 

The group properties of a heat transfer equation with a source
\begin{equation} \label{eq0}
u_t = \left( K(u) u_{x} \right) _{x} + Q(u),
\end{equation}
were considered in [4], and all choices of $K(u)$ and $Q(u)$ which extend the
symmetry group admitted by the general case of equation (1) were identified. 
In this
paper we consider two partial cases of nonlinearities from the
classification in [4] :
\begin{equation} \label{eq1I}
u_t = \left(  u^{\sigma} u_x \right) _x \pm u^{n},
\ \ \sigma ,\  n =const,
\end{equation}
\begin{equation} \label{eq2I}
u_{t} =  u_{xx}  + \delta u \ln u , \ \ \delta = \pm 1,
\end{equation}
together with linear case
\begin{equation} \label{eq3I}
u_t = u_{xx} ,
\end{equation}
whose group properties were known by S.~Lie.
For all cases we construct difference equations
and meshes (lattices) that admit the same Lie group of point
transformations as their continuous limits.
\par
We recall that Lie point symmetries yield a number of useful
properties of differential equations [13,16,10]:
\par
a) A group action transforms the complete set of solutions into itself;
so it is possible to obtain new solutions from a given one.
\par
b) There exists a standard procedure to obtain the whole set of invariants
and differential invariants for a symmetry group; it yields the invariant
representation of the differential equations and the forms of invariant
solutions in which they could be found (symmetry reduction of PDE).
\par
c) For ODEs the known symmetry yields the reduction of the order;
if the dimension of symmetry is equal to (or greater then) the order of ODE,
then we have a complete integrability.
\par
d) The invariance of PDEs is a necessary condition for the application
of Noether's theorem on variational problems to obtain conservation laws.
\par
e) It should be noticed that Lie point transformations have a clear
geometrical interpretation  and one can construct the orbits of a group in a
finite
dimensional space of independent and dependent variables.

\par

The structure of the admitted group
essentially effects the construction of equations
and grids. Group transformations can break the geometric structure of the
mesh that influences the approximation and other
properties of a difference equation. Early contributions to the construction 
of the
difference grids based on the symmetries of the initial difference
model are [6,8]. Classes of transformations that
conserve uniformity, orthogonality, and other properties of meshes 
will be defined below.

\par
In accordance with equation (\ref{eq0}) we consider Lie point
transformations in a space with two independent variables: $t$ and
$x$. Let
\begin{equation} \label{operator}
X = \xi ^{t} {\partial \over \partial t} + \xi ^{x} {\partial \over \partial x}
+\eta {\partial \over \partial u} + \cdots  ,
\end{equation}
\noindent be an operator of a one-parameter transformation group. Dots denote
prolongation of the operator on other variables used in the given
differential equation:
\begin{equation} \label{neq6}
F(x,t,u_t,u,u_x,u_{xx}) = 0 .
\end{equation}
The group generated by (\ref{operator}) transforms a point
$(x,t,u,u_t,u_x,u_{xx})$ to a new one 
$ (x^*, t^*, u^*, u^*_t$, $u^*_x, u^*_{xx})$
together with equation (\ref{neq6}). This situation changes when 
applying 
Lie points transformations to difference equations. Let
\begin{equation} \label{neq7}
F(z) = 0
\end{equation}
be a difference equation defined on some finite set of points
$z^1,z^2,\ldots $ (difference
stencil) on a mesh.
In contrast to the point $(x,t,u,u_t,u_x,u_{xx})$, the "difference point" -
the difference
stencil has its own geometrical structure.
\par
Let
\begin{equation} \label{neq8}
\Omega(z,h) = 0
\end{equation}
be an equation that define a difference stencil and a  mesh. As an example it will  be a uniform mesh if the left
step (spacing) equals the right step:
\begin{equation} \label{neq9}
h^+=h^- .
\end{equation}
\par
The invariance of the difference equation (\ref{neq7}) depends on the
invariance of
(\ref{neq8}), since the latter must be included 
in the general condition of     invariance:
\begin{equation} \label{neq10}
\left\{
\begin{array}{rcl}
XF(z)_{|_{(7)(8)}}& =& 0,\\
X\Omega (z,h)_{|_{(7)(8)}}& =& 0.\\
\end{array}
\right.
\end{equation}
\par
Relations between the two conditions in  (\ref{neq10}) depends on whether 
$ \Omega$, $\xi^t$, $\xi^x$ depend on solution 
or not. If $ \Omega_u=\xi^t_u=\xi^x_u=0$, then the conditions (\ref{neq10})
could be considered independently.
\par
Thus, what makes our approach [5-8] special is 
the inclusion of the second equation of (\ref{neq10}) in the 
conditions of invariance,
admitting the whole set of properties a)--e) stated for
equations (\ref{neq7}) and (\ref{neq8}).
\par
There exist a few ways to avoid transformations of the difference stencil
and, consequently, transformations of a mesh. One way is to restrict
transformations to the case when
independent variables are not changed: $\xi^t=\xi^x=0$, 
yielding any mesh invariant (see [12]). But this restriction is 
very strong and would exclude most symmetries of 
physical problems.
\par
Another approach is connected to evolutionary vector fields. It is known
[10,16] that the symmetry operator    could be represented as
\begin{equation} \label{neq11}
{\bar X} = \left(\eta -  \xi^{t} u_{t}  - \xi^{x}u_{x} \right) {\partial\over
\partial u} + \cdots
\end{equation}
which could be viewed as representative of a factor
algebra by the ideal
\begin{equation}
\xi^t(z) D_t + \xi^x(z) D_x,
\end{equation}
where $ D_t ={\displaystyle \frac{\partial}{\partial t}}+
{\displaystyle u_t \frac{\partial}{\partial u}} + \ldots,\quad
{\displaystyle D_x =\frac{\partial}{\partial x}}+
{\displaystyle u_x\frac{\partial}{\partial u}}\ldots$,
both admitted by all equations (\ref{neq6}). For differential equations
this approach gives an equivalent result but with some losses in points
 a)--e).  In particular we lose the geometric sense even
for point transformations due to the realization of group transformations 
in 
infinite dimensional spaces 
$(t,x,u_t,u_x,u_{tx},u_{xx},u_{tt},\ldots)$, there
is no procedure
for calculating invariants of 
operators of the type given by (\ref{neq11}), etc.
\par
As it was shown in [5],  the exact representation of the operator of total
differentiation $D_x$  in a space of difference variables is given by the
following operators (the same is also true for $D_t$):
\begin{equation} \label{eq13n}
\begin{array}{l}
D^+ ={\displaystyle  {\partial \over \partial x}} +
\tilde {\dd{D}}(u) {\displaystyle {\partial \over \partial u}} + \cdots, 
\qquad 
{\displaystyle \tilde{\dd D}}={\displaystyle \sum_{n=1}^{\infty}
\frac{(-h)^{n-1}}{n} {\dd D}^n} , \\
\\
D^- = {\displaystyle \frac{\partial}{\partial x}} + \tilde{\df D}(u)
{\displaystyle \frac{\partial}{\partial u}} + \cdots, 
\qquad 
 {\displaystyle \tilde{\df D}}={\displaystyle \sum_{n=1}^{\infty}
\frac{h^{n-1}}{n} {\df D}^n} , 
\end{array}
\end{equation}
where 
$\dd D $ and $  \df D$ are right and left 
difference derivative operators on a uniform mesh.
As difference operators are related to corresponding shift operators
$\dd S, \df S$:
$$
     \dg D = \pm \frac{(\dg S -1)}{h},
$$
we obtain another representation for $\tilde{\dg D}$:
\begin{equation} \label{eqnewt}
\tilde{\dd {D}}=\sum_{n=1}^{\infty}\frac{(- h)^{n-1}}{n} {\dd D}^n=
\sum_{n=1}^{\infty}{\frac{(-1)^{n-1}}{nh}(\dd {S}-1)^n},
\end{equation}
\begin{equation}
\tilde{\df {D}}=\sum_{n=1}^{\infty}{\frac{(1-\df{S})^n}{nh}}.
\end{equation}

\par
The group transformation for operators (\ref{eq13n}) can be obtained by
the exponential mapping or by means of the so called Newton series (see [5]).
It is
important to notice that every difference equation on a regular  mesh  admits
the operators $D^+$ and $D^-$, which do not change a mesh. It was shown [5], that
a family of operators $\xi(z) D^{\pm}$ forms an ideal in the Lie algebra of
operators (\ref{operator}),
since it is possible to rewrite these operators as evolutionary vector fields
\begin{equation} \label{eqn12}
\bar X = \left( \eta - {\xi}^t  D^+_{t}(u) - \xi^x D^+_x (u) \right)
{\partial \over \partial u} + \cdots . 
\end{equation}
(We follow here the right semiaxis representation $D^+$. The same formulation
can be done with help of left semiaxis representation $D^-$.)

It is important to notice that the representation (\ref{eq13n}) is true 
only for regular (uniform) meshes. Thus, the evolutionary vector fields 
(\ref{eqn12}) are applicable only for
the groups which do not change uniformity of a mesh. So, one can
not apply them for modern numerical methods with moving meshes, self adaptive
meshes or   multi-grid methods etc.

Let us consider a one-parameter transformation group which is generated 
by the operator 
\begin{equation} \label{eq17n}
X  = \xi \frac{\partial}{\partial x} + 
\eta {\partial \over \partial u} + \cdots .  
\end{equation}
We prolong (\ref{eq17n}) on the right and left steps $h^{+}, h^{-}$ 
by means of relations $h^{+} = x^{+}- x$ and $h^{-} = x - x^{-}$, where
$f^{ \pm} \equiv { \dpm S} (f)$:
\begin{equation} \label{eq17p}
X  = \xi \frac{\partial}{\partial x} + 
\eta {\partial \over \partial u} + \cdots +
[\dd{S}(\xi) - \xi]{\partial \over {\partial h^+}} + 
[\xi - \df{S}(\xi)]\frac{\partial}{\partial h^-}  . 
\end{equation}

   From (\ref{eq17p}) it is easy to obtain the invariance condition for 
a uniform mesh
in a given direction. Let (\ref{neq9}) be invariant with
respect to (\ref{eq17p}), then
\begin{equation} \label{eqn13}
\dd{S}(\xi) - 2\xi +\df{S}(\xi) =0 \quad {\mbox{or}} \quad \dd{D}\df{D}(\xi)=0.
\end{equation}

\par
Condition (\ref{eqn13}) is a strong limitation on the admitted group.
In addition
the coefficients of (\ref{eqn12}) are the  power  series  of
$\dg D$ or $\dg S$, since one should consider the whole set of 
mesh points and not a stencil only.
\par
Let us illustrate the above with a simple example. The ODE
\begin{equation} \label{eqn14}
u_{xx}=u^2,
\end{equation}
can be viewed as the stationary case of equation (\ref{eq1I})
with $\sigma=0$, $n=2$. Equation
(\ref{eqn14}) has the following Lie point
symmetries:
\begin{equation} \label{eqn15}
X_1=\frac{\partial}{\partial x}; \quad
X_2=x\frac{\partial}{\partial x} -2u \frac{\partial}{\partial u} .
\end{equation}
As a difference analog of (\ref{eqn14}) we consider
\begin{equation} \label{eqn16}
\frac{u^+-2u+u^-}{h^2} = u^2.
\end{equation}
on a uniform mesh
\begin{equation} \label{eqn17}
h^+=h^-,
\end{equation}
where $u^+=\dd{S}(u), \quad u^-=\df{S}(u)$.

The equations (\ref{eqn16}) and (\ref{eqn17}) use a three-point
stencil or subspace $(x, x^+, x^-, u, u^+, u^-)$ and 
the operators (\ref{eqn15}) have the following prolongation for the shifted
points
of the difference stencil:
\begin{equation} \label{eqn23}
\begin{array}{c}
{  \displaystyle 
X_1=\frac{\partial}{\partial x}+ \frac{\partial}{\partial {x^+}}+
\frac{\partial}{\partial {x^-}}; 
\qquad 
X_2= x\frac{\partial}{\partial x} + x^+
\frac{\partial}{\partial {x^+}}  +
 x^- \frac{{\partial}}{{\partial} {x^-}}  } \\
\\
{  \displaystyle 
 -2u\frac{\partial}{{\partial}{u}} -
2u^+ \frac{\partial}{\partial {u^+}} -2u^-\frac{\partial}{\partial {u^-}} +
 h^+{\partial \over {\partial h^+}} +
h^- \frac{\partial}{\partial h^-}. } 
\end{array}
\end{equation}

The symmetry algebra (\ref{eqn15}) acts on the space $(t,x)$,
so the coefficients
of (\ref{eqn23}) have the same form in different points of the stencil. The
prolongation
forms for $h^+$ and $h^-$ are easily derived from the relations
$ h^+ = x^+ - x,\quad h^-= x^- - x$.
\par
It is easy to verify that the  equations (\ref{eqn16}) and (\ref{eqn17}) are
invariant under
the action of (\ref{eqn23}):
\begin{equation} \label{eqn18}
\begin{array}{rcl}
X_2 \left({\displaystyle \frac{u^+ -2u + u^-}{h^2}}
 - u^2 \right)_{|_{(\ref{eqn16})}}& =&0, \\
\\
X_2 \left( h^+ -h^- \right)_{|_{(\ref{eqn17})}}& =&0,
\end{array}
\end{equation}
(Operator $X_1$  leaves (\ref{eqn16}) and (\ref{eqn17}) unchanged).
\par
It follows that the system (\ref{eqn16}) and (\ref{eqn17}) 
has the same Lie point symmetry as its continuous limit.
Notice that the invariance conditions (\ref{eqn18}) are mutually independent.
\par
As for the continuous case one can easily  calculate the finite difference
invariants for (\ref{eqn23}) by solving the standard system:
\begin{equation} \label{eqn19}
X_i \  I (x,h^+,h^-,u,u^+,u^-)=0,\qquad i=1,2. 
\end{equation}
The solution of (\ref{eqn19}) yields the whole set of difference invariants
\begin{equation} \label{eqn20}
I^1=\frac{h^+}{h^-},\quad
I^2=\frac{u^+}{u},\quad
I^3=\frac{u^-}{u},\quad
I^4=(h^+)^2 u;
\end{equation}
It follows that the difference model (\ref{eqn16}) and (\ref{eqn17}) 
can be represented by
means of the invariants (\ref{eqn20}) as $I^2 + I^3 - 2 = I^4$ and $I^1 = 1$.

\par
Let us now consider the evolutionary vector fields for the difference
equation  (\ref{eqn16}) 
(we consider the right semiaxis representation $D^+$). 
We prolong the operator (\ref{eq13n}) on all points of a given stencil 
$(x,x^{+},x^{-},u,u^{+},u^{-})$:
\begin{equation} \label{eqn21n}
D^{+} = 
{\ddx} + 
{\partial \over \partial x^{+} } + 
{\partial \over \partial x^{-} } + 
u_{x} {\ddu }  + 
u_{x}^{+}  {\partial \over \partial u^{+} } + 
u_{x}^{-}  {\partial \over \partial u^{-} } ,
\end{equation}
where 
$$
u_{x} \equiv 
{\displaystyle \sum_{n=1}^{\infty}\frac{(-h)^{n-1}}{n}{\dd D}^n (u) }, 
\quad 
u_{x}^{+}  \equiv 
{\displaystyle \sum_{n=1}^{\infty}\frac{(-h)^{n-1}}{n}{\dd D}^n (u^{+}) } , 
\quad 
u_{x}^{-} \equiv 
{\displaystyle \sum_{n=1}^{\infty}\frac{(-h)^{n-1}}{n}{\dd D}^n (u^{-}) } .
$$  
The evolution vector fields (\ref{eqn12}) 
will have the following forms:
\begin{equation} \label{eqn21p}
\begin{array}{lcl}
  { \dh  {\bar{X}} { }_1   } &   = &  - X_{1} + D^{+}  = 
{\displaystyle 
u_{x} {\ddu} + 
u_{x}^{+}  {\partial \over \partial u^{+} } + 
u_{x}^{-}  {\partial \over \partial u^{-} }  }; \\
  \\
 {  \dh   {\bar{X}} { }_2 }  &  =  & - X_{2} + x D^{+} = 
{\displaystyle 
- h^{+} {\partial \over \partial h^{+} } - 
h^{-} {\partial \over \partial h^{-} }  } \\
\\
& & 
{\displaystyle 
+ (2u + x u_{x}) {\ddu} + 
(2u^{+}  + x u_{x}^{+}) {\partial \over \partial u^{+} }+ 
(2u^{-}  + x u_{x}^{-}) {\partial \over \partial u^{-} } }. 
\end{array}
\end{equation}

%

\par
It is  not  easy  to  check  the  invariance  conditions  of  equation
(\ref{eqn16})  for the operators
(\ref{eqn21p}) because one should use not only the equation (\ref{eqn16}),
but all its sequences
obtained by shifting to the right. A harder question is how to
produce the finite difference invariants (\ref{eqn20}) by means of
({\ref{eqn21p}). 
That is why we prefer the first
classical way for Lie point symmetries and include a mesh in the invariance
condition ({\ref{neq10}) (a very similar approach for the semi-discretized
nonlinear
heat equation was introduced recently  in [3]). An additional convincing reason
to apply the classical representation of Lie point symmetries springs from
comparison of unsuitable evolutionary vector fields approaches to invariant
variational problems, developed in [5], and a clear classical way to construct
Noether type theorems for difference equations [7].
\par
Another approach to the symmetry of difference equations 
on a fixed uniform mesh was introduced 
in [9]. However, that way is only applicable to linear equations; 
moreover, it requires
knowing a complete set of solutions of difference equations. The newly
introduced
approach [11] deals with evolutionary vector fields on a uniform mesh. The
advantage of the last two approaches seems to be in the potential of finding
non-point symmetries of difference equations which are not available in 
their continuous limits.

\par
Returning to equation (\ref{eq0}), we notice that a transformation defined by
(\ref{operator})
conserves uniformity of a grid in $t$ and $x$ directions, if
\begin{equation} \label{ct}
{ \da{D} }{ \db{D}  }( \xi ^{t} ) = 0 ,
\end{equation}
\begin{equation} \label{ch}
{ \dd{D} }{ \df{D} }  ( \xi ^{x} ) = 0,
\end{equation}
where $\da D$ and $\db D$ are right and left difference operators in $t$ direction.


\par
Conditions (\ref{ct}) and (\ref{ch}) are not sufficient to describe 
the invariance 
of an orthogonal mesh. For an orthogonal mesh to conserve its orthogonality
under the
transformation, it is necessary and sufficient that [6,8]:
\begin{equation} \label{cht}
{\dd{D} }( \xi ^{t} ) = -{ \da{D}} ( \xi ^{x} ) .
\end{equation}
\par
When condition (\ref{cht}) is not satisfied for a given group, the
flatness of the layer of a grid in some direction is rather important.
For evolution equations it is significant to have flat time layers, since
otherwise, after transformations, some domains of a space could be
in the future, while others in the past. We have a simple
criterion [6,8] for the preserving the flatness of the layer of a grid
in the time direction under the action of a given operator
(\ref{operator}):
\begin{equation} \label{cht2}
{\dd{D}} { \da{D}} ( \xi ^{t} ) = 0.
\end{equation}

\par

So, the conditions (\ref{ct})--(\ref{cht2}) provide a geometry of grids
that is based on the Lie group symmetry. These conditions will be used
in what follows.

\bigskip

\noindent{\bf 2. Invariant model for the equation $u_{t} =
(u^{\sigma}u_{x})_{x} \pm u^{n}$  }

\medskip
\noindent The equation
\begin{equation} \label{eq33}
u_{t}=(u^{\sigma}u_{x})_{x}  \pm u^{n},\quad \sigma,n=const,
\end{equation}
admits a 3-parameter symmetry group. This group can be represented by the
following infinitesimal operators [4]:
\begin{equation} \label{op33}
X_{1}= {\ddt};\quad X_{2}={\ddx};\quad X_{3}= 2(n-1)t{\ddt}+(n-\sigma
-1)x{\ddx} -2u{\ddu} .
\end{equation}
The set (\ref{op33}) satisfies all conditions (\ref{ct})--(\ref{cht2}).
So, we can use an orthogonal grid that is uniform in the $t$ and $x$
directions.
Let us consider the set of operators (\ref{op33}) in the space
$( t,  \hat{t},  x,  h^{+},  h^{-},  u,  u_{+}, 
u_{-},   \hat{u},  \hat{u}_{+},  \hat{u}_{-})$ that
corresponds to the stencil shown in Fig. 1.

\begin{figure}[h!]  \label{grort}


\begin{picture}(300,100)




%

\put(105,0){\begin{picture}(200,100)%

\put(30,30){\line(1,0){140}}

\put(30,70){\line(1,0){140}}

\put(100,30){\line(0,1){40}}

\put(30,30){\circle*{5}}

\put(30,70){\circle*{5}}

\put(170,70){\circle*{5}}

\put(170,30){\circle*{5}}

\put(100,30){\circle*{5}}

\put(100,70){\circle*{5}}

\put(5,75){$(x-h,\hat{t},\hat{u}_{-})$}
\put(85,75){$(x,\hat{t},\hat{u})$}
\put(145,75){$(x+h,\hat{t},\hat{u}_{+})$}
\put(5,20){$(x-h,t,u_{-})$}
\put(85,20){$(x,t,u)$}
\put(145,20){$(x+h,t,u_{+})$}













%

\end{picture}}

%


%

\end{picture}

\caption{The stencil of the orthogonal mesh.}
\end{figure}
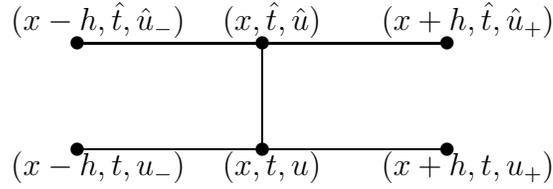

There are 7 difference invariants of the Lie algebra (\ref{op33}):
\begin{equation} \label{in33}
{\tau ^{ n -\sigma -1 \over 2(n-1) } \over h };
\quad \tau u^{n-1} ; \quad { \hat{u} \over u} ;
\quad {u_{+} \over u}; \quad { u_{-} \over u} ;
\quad {\hat{u}_{+} \over \hat{u}}; \quad { \hat{u}_{-} \over \hat{u}} .
\end{equation}
The small number of symmetry operators (\ref{op33}) provides us with a large
number of
difference invariants (\ref{in33}). Thus we are left with some additional 
degrees of freedom in invariant difference 
modeling of (\ref{eq33}). By means of the invariants (\ref{in33}), we could
write
the following explicit scheme for (\ref{eq33}):
\begin{equation} \label{sh33}
{ \hat{u} - u \over \tau } =
{ 1  \over h} \left(  \left({u_{+} +u \over 2} \right)^{\sigma} \dh  u { }_{x}
- \left({u+ u_{-}\over 2 }\right)^{\sigma} \dh u { }_{\bar x}  \right)  \pm
  u^{n},
\end{equation}
where $\dh u { } _{x}  = {\displaystyle {u_{+} - u \over h }}$, $\dh u { }_{ \bar x}  =
{\displaystyle {u - u_{-} \over h}}$.
\medskip

\bigskip
This scheme is not unique and one could construct another form of invariant
difference equation. For example an implicit scheme could be as follows:
\begin{equation} \label{shsh}
{ \hat{u} - u  \over \tau } =
{ 1\over h^{2}}
( \hat{u}_{+}^{\sigma +1 } -
2 \hat{u}^{\sigma +1 } +
 \hat{u}_{-}^{\sigma +1 } ) +
 \hat{u}^{n}
\end{equation}
Notice, that continuous limit of the last difference equation
\begin{equation} \label{eq33a}
u_{t} =
( {u}^{\sigma +1 } )_{xx} +  u^{n}.
\end{equation}
is equivalent to the equation (\ref{eq33})
up to the scaling of $x$. But scheme (\ref{shsh}) is not equivalent to the
scheme
(\ref{sh33}), because there is no point transformation that relates them.
In [18] Samarskii { \it et all}  considered
the case $n = \sigma + 1$, $\sigma >0$ and found a finite-difference blow-up
solution for the equation
(\ref{shsh}) that is invariant with respect to the operator
\begin{equation} \label{operat}
X^*_3 = ( t - T  ) { \partial \over \partial t } -
{ 1 \over \sigma  } u { \partial \over \partial u },
\end{equation}
where $T$ is constant. The operator (\ref{operat}) defines
a subgroup which
is equivalent to the self-similar subgroup with the operator
$X_{3}$ of the set (\ref{op33}). Let us find the solution of the
problem (\ref{shsh}) in the invariant form:
$$
u = \left( 1 - { t \over T  } \right) ^{  \displaystyle  - { 1 \over
\sigma }}  \theta (x).
$$
This solution is sought ([18]) on the time mesh with 
the infinite number of decreasing steps 
\begin{equation} \label{m1}
\tau _{j} =   \sigma   { \rho - 1  \over \rho ^{\sigma +1} } 
\rho ^{- \sigma j }, \qquad j = 0,1,2,...,
\end{equation}
where $\rho > 1$ is constant. This mesh corresponds to the blow-up time 
$$
 T  = {  \sigma \over \rho }  { \rho - 1  \over \rho ^{\sigma } -1  } . 
$$ 
For the function $\theta(x)$ we have the equation
\begin{equation} \label{m3}
( \theta ^{\sigma + 1} ) _{ \bar{x} x } +
 \theta ^{\sigma + 1}  =
{ 1 \over \sigma  } \theta .
\end{equation}
The solution of the problem (\ref{m3})
was found in [18] for the case $\sigma = 2$. Let us fix an arbitrary $M > 0$
and let $h = 2 \sin \left( { \displaystyle { 3 \pi \over 2(M+1) } } \right)$.
In this case the localization length equals
$$
l_{h} = { 3 \pi h \over 2 } \left( \arcsin {h\over 2 } \right) ^{-1} , \ \ 0 <
h \leq 2
$$
(see [18]). Then, one can verify that the solution of the problem
(\ref{m3}) in the points $x_k=kh$ has the form
\begin{equation} \label{s2}
\theta _{kh} =  \sqrt{2} \left(  3 \left( 1 - { 4\over h^{2}} \sin^{2} {a_{h}
h\over 2} \right) \right) ^{ -{1/2} }
\sin ( a_{h} kh ),\ \
k= 0, 1, ...,M+1,
\end{equation}
where $a_{h} = \pi / l_{h} $.

The obtained function $u$ gives the blow-up
solution of the problem in the case $\sigma =  2$, $n = 3$, $l= l_{h}$.
This solution tends to infinity in all points of the space grid,
conserving the structure.
As $h \rightarrow 0$, the difference solution (\ref{s2}) tends
to the solution of the ordinary differential equation:
$$
\theta (x) = \left( { 3 \over 4 } \right) ^{1/2} 
\sin \left( { x \over 3 } \right) , \qquad  0< x < l_{0} = 3 \pi .
$$

\bigskip

\noindent{
\bf
3. Invariant difference model for the semilinear heat transfer equation }

\medskip
\noindent The semilinear heat transfer
equation
\begin{equation} \label{eq53}
u_{t}=u_{xx} + \delta  u \ln u, \ \  \delta = \pm 1,
\end{equation}
admits the 4-parameter Lie symmetry group of point
transformations [4], that corresponds to the following
infinitesimal operators:
\begin{equation} \label{op53}
X_{1}= {\ddt}; \quad X_{2} = {\ddx}; \quad  X_{3} = 2 {\Large e}^{\displaystyle
 \delta t } {\ddx} - \delta {\Large e}^{\displaystyle  \delta t} x u{\ddu};
\quad X_{4} = {\Large e}^{\displaystyle  \delta t } u {\ddu} .
\end{equation}
Before constructing a difference equation and a grid that
approximate (\ref{eq53}) and inherit the whole Lie algebra (\ref{op53}),
we should first check condition (\ref{cht}) for
the invariance of orthogonality. The operators $X_{1}$, $X_{2}$
and $X_{4}$ conserve the orthogonality, but $X_{3}$ does not: condition
(\ref{cht}) is not true for the last operator. Consequently an orthogonal 
mesh cannot be used for the invariant modeling of (\ref{eq53}).

The conditions (\ref{cht2}) are true for the complete set (\ref{op53}), so 
it is possible to use a nonorthogonal grid with flat time layers. An example 
of a grid with with flat time layer is shown in Fig. 2.

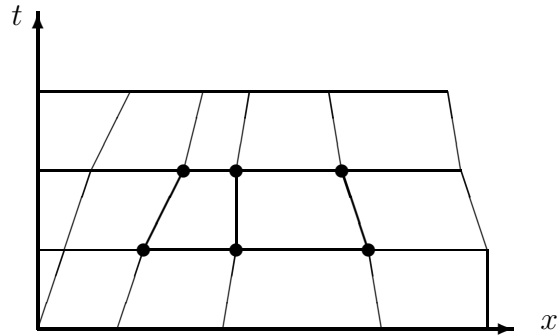
\begin{figure}[h!]



\begin{picture}(300,160)





%

\put(120,0){\begin{picture}(200,150)%
{ \thicklines

\put(10,20){\vector(1,0){180}}

\put(10,20){\vector(0,1){120}}

}

\put(10,50){\line(1,0){170}}

\put(10,80){\line(1,0){160}}

\put(10,110){\line(1,0){155}}

\put(10,20){\line(1,3){10}}

\put(40,20){\line(1,3){10}}

\put(80,20){\line(1,6){5}}

\put(140,20){\line(-1,6){5}}

\put(180,20){\line(0,1){30}}

\put(20,50){\line(1,3){10}}

{\thicklines

\put(50,50){\line(1,2){15}}

\put(85,50){\line(0,1){30}}

\put(135,50){\line(-1,3){10}}

\put(50,50){\line(1,0){85}}

\put(65,80){\line(1,0){60}}

}

\put(180,50){\line(-1,3){10}}

\put(30,80){\line(1,2){15}}

\put(65,80){\line(1,4){7,5}}

\put(85,80){\line(1,6){5}}

\put(125,80){\line(-1,6){5}}

\put(170,80){\line(-1,6){5}}

\put(50,50){\circle*{5}}

\put(85,50){\circle*{5}}

\put(135,50){\circle*{5}}

\put(65,80){\circle*{5}}

\put(85,80){\circle*{5}}

\put(125,80){\circle*{5}}

\put(0,135){$t$}

\put(200,20){$x$}

\end{picture}}
%
%
\end{picture}

\caption{An evolutionary mesh with flat time-layers.}
\end{figure}

A possible reformulation of equation (\ref{eq53}) by using the four 
differential invariants in the subspace 
$(t, x, u, u_{x}, u_{xx}, dt, dx, du)$:
$$
J^{1} = dt; \quad     J^{2} = \left( { u_{x} \over u } \right) ^{2} - { u_{xx}
\over u } ;
\quad J^{3} = 2 { u_{x} \over u } + { dx \over dt } ; 
\quad    
J^{4} = {du \over u dt } - \delta  \ln u + { 1\over 4} \left( {dx \over dt}
\right)^{2} ;
$$
is given by the system:
$$
\left\{
\begin{array}{rcl}
J^{3}& =& 0 ; \\
J^{4}& =& J^{2} ;\\
\end{array}
\right.
$$
or
\begin{equation} \label{sys53}
\left\{
\begin{array}{l}
{\displaystyle {dx\over dt}  = -2 {u_{x} \over u }    } ;\\
\\
{\displaystyle {du\over dt}  = u_{xx} + \delta u \ln u - 2 { u_{x}^{2} \over u
}  }.\\
\end{array}
\right.
\end{equation}
So, the structure of the group (\ref{op53}) suggests the use of two
evolution equations.
\par
As the next step, we will find difference invariants for the set
$X_{1}$--$X_{4}$ of the group (\ref{op53}). These invariants are necessary
for the approximation of the system (\ref{sys53}). We will use
the 6-point difference stencil, as shown in Fig. 3

\begin{figure}[h!]   \label{grnonor}


\begin{picture}(300,100)





%

\put(110,0){\begin{picture}(200,100)%
\put(20,30){\line(1,0){160}}

\put(30,70){\line(1,0){140}}

\put(80,30){\line(1,2){20}}

\put(20,30){\circle*{5}}

\put(30,70){\circle*{5}}

\put(170,70){\circle*{5}}

\put(180,30){\circle*{5}}

\put(80,30){\circle*{5}}

\put(100,70){\circle*{5}}

\put(0,75){$(\hat{x}-\hat{h}^{-},\hat{t},\hat{u}_{-})$}
\put(85,75){$(\hat{x},\hat{t},\hat{u})$}
\put(140,75){$(\hat{x}+\hat{h}^{+},\hat{t},\hat{u}_{+})$}
\put(65,20){$(x,t,u)$}
\put(150,20){$(x+h^{+},t,u_{+})$}
%






%

\end{picture}}

%


\put(90,20){$(x-h^{-},t,u_{-})$}

\end{picture}

\caption{The stencil of the evolutionary mesh.}
\end{figure}
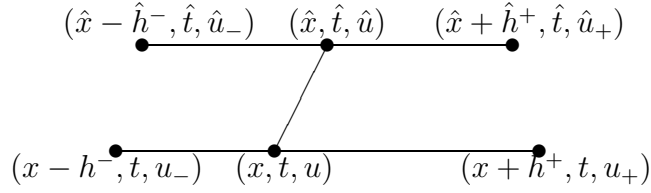

\noindent on which we will approximate the system (\ref{sys53}). The stencil 
defines the difference subspace $( t, \hat{t}, x, \hat{x}, h^{+}, h^{-},
\hat{h}^{+},
\hat{h}^{-}, u,$$ u_{+},$$ u_{-},\hat{u},$$ \hat{u}_{+},$$ \hat{u}_{-} )$.
The group (\ref{op53}) has the following difference invariants in this
subspace:
$$
I^{1} =  \tau; \quad I^{2} =  {h}^{+}; \quad I^{3} = {h}^{-};
\quad I^{4} = \hat{h}^{+}; \quad I^{5} = \hat{h}^{-};
$$
$$
I^{6} = ( \ln u) _{x}  - ( \ln u )_{\bar{x}} ;
\qquad I^{7} = ( \ln \hat{u} ) _{x}  - ( \ln \hat{u} )_{\bar{x}} ;
$$
$$
I^{8} = \delta \Delta x + 2 ({\Large e}^{\displaystyle \delta \tau} - 1)
\left( { h^{-} \over h^{+} + h^{-} } ( \ln u )_{x}  +
    { h^{+} \over h^{+} + h^{-} } (  \ln u )_{\bar{x}}  \right) ;
$$
$$
I^{9} = \delta \Delta x + 2 (1 - {\Large e}^{\displaystyle - \delta \tau} )
\left( { \hat{h}^{-} \over \hat{h}^{+} + \hat{h}^{-} } ( \ln \hat{u} )_{x} +
   {\hat{h}^{+} \over \hat{h}^{+} + \hat{h}^{-} } ( \ln \hat{u} )_{\bar{x}}
\right) ;
$$
$$
I^{10} = \delta (\Delta x) ^{2} + 4 (1 - {\Large e}^{\displaystyle - \delta
\tau} )
\left( \ln \hat{u} - {\Large e}^{\displaystyle  \delta \tau} \ln u \right) ;
$$
where ${\displaystyle \Delta x = \hat{x}-x }$,
${\displaystyle  (\ln u)_{x}={ \ln u_{+} - \ln u \over h^{+}} }$ ,
${ \displaystyle ( \ln u ) _{\bar x} = { \ln u - \ln u_{-} \over h^{-}} }$ .

\medskip

\par
\noindent To obtain an invariant difference model, it is natural 
to use the difference invariants. An explicit model is given by
$$
\left\{
\begin{array}{l}
{ \displaystyle I ^{8} = 0 }  ;\\
\\
{ \displaystyle I^{10} =  {8 \over \delta} 
{  ({\Large e}^{\displaystyle \delta I^{1} } - 1)^{2} \over
I^{2} + I^{3} } I^{6} } ; \\
\end{array}
\right.
$$
or
\begin{equation} \label{sh53}
\left\{
\begin{array}{l}
{ \displaystyle {\delta \Delta x + 2 ({\Large e}^{\displaystyle \delta \tau} -
1)
\left( { h^{-} \over h^{+} + h^{-} } ( \ln u) _{x}  +
    { h^{+} \over h^{+} + h^{-} } ( \ln u )_{\bar{x}}  \right) = 0 } };\\
\\
{ \displaystyle {\delta (\Delta x) ^{2} + 4 (1 - {\Large e}^{\displaystyle -
\delta \tau} )
\left( \ln \hat{u} - {\Large e}^{\displaystyle  \delta \tau} \ln u \right) =
{8 \over \delta} {  ({\Large e}^{\displaystyle \delta \tau} - 1)^{2} \over
h^{+} + h^{-} }
           \left[ (\ln u) _{x} - ( \ln u)_{\bar{x}} \right] } } .
\end{array}
\right.
\end{equation}

\medskip

\par
As in the continuous case, there is a reduction 
in the difference case. 
When we consider an invariant solution, we have the reduction of the
equation-grid
system to a system of ordinary difference equations. 
One being the difference model for the considered equation, the other 
for the evolution of the grid.

  Let us find the solution of the difference model (\ref{sh53}) which is
invariant with respect to the operator
\begin{equation} \label{oper53}
2 \alpha X_{2} + X_{3}, \ \ \alpha = const.
\end{equation}
${ \displaystyle { u \exp \left(  { \delta {\Large e}^{\displaystyle  \delta t
}
\over
\alpha + {\Large e}^{\displaystyle  \delta t } }
{ x^{2}  \over 4 } \right)} } $,
$\left( {\displaystyle {\Delta x\over
{\Large e}^{\displaystyle  \delta t }
 ( {\Large e}^{\displaystyle  \delta \tau  } -1 )  } -
{ x \over
 \alpha + {\Large e}^{\displaystyle  \delta t } } } \right) $
and $t$ are all the invariants with respect to (\ref{oper53}). 
Therefore
we will seek an invariant solution in the form:
$$
\left\{
\begin{array}{l}
{ \displaystyle { u(x,t)\ \  =\ \   \exp \left( -
 { \delta {\Large e}^{\displaystyle  \delta t } \over
\alpha + {\Large e}^{\displaystyle  \delta t } }
{ x^{2} \ \over 4 } \right)  {\Large e}^{\displaystyle  f(t)  } } } ;\\
\\
{\displaystyle {\Delta x\over
{\Large e}^{\displaystyle  \delta t }
 ( {\Large e}^{\displaystyle  \delta \tau  } -1 )  } \ \ = \ \
{ x \over
 \alpha + {\Large e}^{\displaystyle  \delta t } } + g(t)  } .\\
\end{array}
\right.
$$
\noindent Substituting this form of the solution into 
(\ref{sh53}), we obtain a system of ordinary difference equations to
determine $f(t)$ and $g(t)$:
$$
\left\{
\begin{array}{l}
{ \displaystyle { f(t+\tau) - {\Large e}^{\displaystyle  \delta \tau } f(t)
\over
{\Large e}^{\displaystyle  \delta \tau  }
 ( {\Large e}^{\displaystyle  \delta \tau  } -1 )  }  } \ \ = \ \
{\displaystyle {  - {1\over 2}
{  {\Large e}^{\displaystyle  \delta t }  \over
\alpha + {\Large e}^{\displaystyle  \delta t  } } } };\\
\\
g(t)\ \  =\ \  0.\\
\end{array}
\right.
$$
\par
\noindent The solution of this system yields the solution of the difference
equation (\ref{sh53}):
$$
u(x,t) = \exp \left(   {\Large e}^{\displaystyle  \delta t  }
\left( f(0) - {   {\Large e}^{\displaystyle  \delta \tau  }  - 1 \over 2 }
\sum_{j=1}^{n-1}   { {\Large e}^{\displaystyle -  \delta t_{j} } \over
 1+ \alpha  {\Large e}^{\displaystyle -  \delta t_{j}   } } \right) -
 { \delta {\Large e}^{\displaystyle  \delta t } \over
\alpha + {\Large e}^{\displaystyle  \delta t } }
{ x^{2} \ \over 4 } \right)   .
$$
\noindent and the grid
$$
x = x^{0}  { {\Large e}^{\displaystyle  \delta t  } + \alpha
\over 1 + \alpha }.
$$
\noindent Here $x = x^{j}_{i} = x_{i}(t_{j})$ and $ t = t_{j}$.
For $t=0$ the grid can be arbitrary, but if a regular grid is used 
the grid will be regular on every time layer.

\par
The obtained solution is the solution of the Cauchy problem with
initial conditions:
$$
u(x,0) =    \exp \left( f(0) -
 { \delta {\Large e}^{\displaystyle  \delta t } \over
\alpha + {\Large e}^{\displaystyle  \delta t } }
{ x^{2} \ \over 4 } \right)   .
$$

\bigskip

\noindent {\bf 4. Invariant discrete version of the linear heat equation }

\medskip
\noindent The linear heat transfer equation
\begin{equation} \label{eq54}
u_{t} = u_{xx}
\end{equation}
admits a 6-parameter Lie symmetry group of point transformations, 
corresponding to the following infinitesimal operators:
\begin{equation} \label{op54}
\begin{array}{c} 
{ \displaystyle 
X_{1} = {{\partial \over \partial t}} ;\qquad X_{2} = {{\partial \over \partial
x}} ; \qquad X_{3} = {2t} {{\partial \over \partial x}} {-} xu {{\partial \over
\partial u}} ; } \\
\\
{ \displaystyle 
X_{4} = {2t} {{\partial \over \partial t}} + {x} {{\partial \over \partial x}}
; \quad X_{5} = {4t^{2}} {{\partial \over \partial t}} + {4tx} {{\partial \over
\partial x}} - ( {x^{2}} + {2t} {)u} {{\partial \over \partial u}} ;
\quad X_{6} = {u} {{\partial \over \partial u}} ; } 
\end{array} 
\end{equation}
\noindent and an infinite-dimensional symmetry obtained from 
the linearity of the equation (\ref{eq54}):
$$
X^{*} = {a(x,t)} \ {{\partial \over \partial u}} ,
$$
where $a(t,x)$ in an arbitrary solution of the equation (\ref{eq54}).

Now we are in a position to show 
that the invariant difference model for the linear
heat
transfer equation cannot be constructed on an orthogonal grid. The model
\begin{equation}  \label{eq55a}
{ \hat{u} - u \over \tau } = { u_{+} -2u + u_{-} \over  h
^{2} }
\end{equation}
on the regular orthogonal mesh, which is used as an 
invariant model in [2], (see also
[1], p.~363), 
actually does not admit operators $X_{3}$ and $X_{5}$.

Let us check for example the symmetry that is described by $X_{3}$.
The operator $X_{3}$ generates the following transformations
\begin{equation}  \label{tr}
\begin{array}{l}
x^{*} = x + 2t \alpha, \\
t^{*} = t , \\
u^{*} = u {\Large e }^{\displaystyle - x \alpha - t \alpha ^{2} }.
\end{array}
\end{equation}
This transformation destroys the orthogonality of the mesh as 
shown in Figure 4.

\begin{figure}[h!]

\begin{picture}(300,120)
\put(60,10){\begin{picture}(130,100)%

\put(10,10){\vector(1,0){120}}
\put(10,10){\vector(0,1){90}}
\put(10,40){\line(1,0){100}}
\put(10,70){\line(1,0){100}}
\put(40,10){\line(0,1){75}}
\put(70,10){\line(0,1){75}}
\put(100,10){\line(0,1){75}}

\put(15,100){$t$}
\put(120,15){$x$}

\put(13,30){$x-h$}
\put(75,30){$x$}
\put(105,30){$x+h$}
\put(75,75){$ \hat{x}$}

\put(40,40){\circle*{5}}
\put(70,40){\circle*{5}}
\put(100,40){\circle*{5}}
\put(70,70){\circle*{5}}
%
%
\end{picture} }
\put(225,10){\begin{picture}(130,100)%
\put(10,10){\vector(1,0){120}}
\put(10,10){\vector(0,1){90}}
\put(10,40){\line(1,0){110}}
\put(10,70){\line(1,0){110}}
\put(25,10){\line(1,2){37}}
\put(55,10){\line(1,2){37}}
\put(85,10){\line(1,2){37}}
\put(15,100){$t$}
\put(120,15){$x$}
\put(40,40){\circle*{5}}
\put(70,40){\circle*{5}}
\put(100,40){\circle*{5}}
\put(85,70){\circle*{5}}
\put(15,30){$x^{*} - h^{*}$}
\put(73,30){$x^{*}$}
\put(103,30){$x^{*} + h^{*}$}
\put(93,75){$ \hat{x}^{*}$}
\multiput(85,70)(0,-9){4}{\line(0,-1){6}}
%
%
\end{picture} }
\end{picture}
\caption{Deformation of the orthogonal mesh.}
\end{figure}
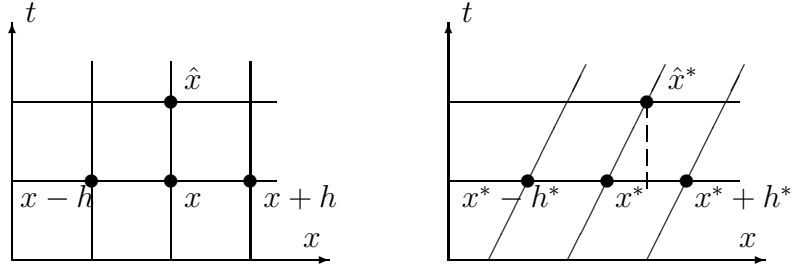

The transformation (\ref{tr}) transforms 
the finite-difference equation (\ref{eq55a}) 
into the following:
\begin{equation} \label{eqor}
{ \hat{u} {\Large e} ^{\displaystyle  - \tau \alpha^{2} }
- u \over \tau } =
{ u_{+} {\Large e}^{ \displaystyle - h \alpha } - 2u +
u_{-} {\Large e } ^{\displaystyle h \alpha }   \over h^{2} },
\end{equation}
which explicitly depends on a group parameter $\alpha$. 
The first differential approximation of the equation (\ref{eqor})
$$
u_{t} = u_{xx} - 4 \alpha u_{x}     + 2 u \alpha ^{2} + O(\tau +
h )
$$
shows explicitly the absence of invariance for the equation
(\ref{eq55a}) on an orthogonal mesh.

Consequently, we have to construct a difference model for (\ref{eq54}) on a
moving mesh. 
With help of the differential invariants in the space 
$( t, x, u, u_{x}, u_{xx}, dt, dx, du)$:
$$
J^{1} = { { dx + 2 {\displaystyle { u_{x} \over {u} } } dt } \over { dt^{1/2}}}
;
\qquad J^{2} = {{du\over u}} + {{1\over 4}} {{dx\over dt}}^{2} + \left( -
{{u_{xx}\over u}} + {u_{x}^{2}\over u^{2}} \right) dt ;
$$
\noindent we can represent the heat equation (\ref{eq54})
as the system:
$$
\left\{
\begin{array}{rcl}
J^{1} &=& 0;  \\
J^{2} &=& 0;  \\
\end{array}
\right.
$$
\noindent or
\begin{equation} \label{sys54}
\left\{
\begin{array}{rcl}
{\displaystyle {dx\over dt} }&=&   -2 {\displaystyle {u_{x}\over u } }; \\
\\
{\displaystyle {du\over dt} }&=& u_{xx} - {2}{\displaystyle {u_{x}^{2}\over u}
}. \\
\end{array}
\right.
\end{equation}
The system (\ref{sys54}) inherits the set of operators $X_{1}$--$X_{6}$,
$X^{*}$.
Implying that it entirely inherits the Lie
symmetry group admitted by the linear heat  equation (\ref{eq54}).
\par
For the difference modeling of the system (\ref{sys54}) we need the whole set
of difference invariants of the Lie symmetry group (\ref{op54}) in the
difference space, corresponding to the chosen stencil $( t, \hat{t}, x,
\hat{x}, h^{+}, h^{-},
\hat{h}^{+}, \hat{h}^{-},u, \hat{u}, u_{+}, u_{-}, 
\hat{u}_{+}, \hat{u}_{-} )$:
$$
I^{1} = { {h^{+}} \over {h^{-}} } ;
\quad I^{2} = { {\hat h^{+}} \over {\hat h^{-}}} ;
\quad I^{3} = { { {\hat h^{+}} { h^{+}}} \over \tau } ;
\quad I^{4} = { \tau^{1/2} \over h^{+} }  {{\hat{u}\over u}} {\exp } \left(
{{1\over 4}} {{(\Delta x)^{2}\over \tau}} \right) ;
$$
$$
I^{5} = {{1\over 4}} {{h^{+2}\over \tau}} - {{h^{+2}\over h^{+} + h^{-}}}
\left( {{1\over h^{+}}} {\ln }  {{u_{+}\over u}}  + {{1\over h^{-}}} {\ln }
 {{u_{-}\over u}} \right) ;
$$
$$
I^{6} = {{1\over 4}} {{\hat{h}^{+2} \over \tau}} + {\hat{h}^{+2} \over
\hat{h}^{+} + \hat{h}^{-}} \left( {1\over \hat{h}^{+}} {\ln } { \hat u_{+}
\over \hat u} + {1\over \hat{h}^{-}} {\ln } {\hat u_{-} \over \hat u} \right) ;
$$
$$
I^{7} = {\Delta xh^{+} \over \tau}  + {2h^{+} \over h^{+} + h^{-}}  \left(
{h^{-} \over h^{+}} {\ln }  {u_{+} \over u}  - {h^{+} \over h^{-}} {\ln }
 {u_{-} \over u}  \right) ;
$$
$$
I^{8} = {{\Delta x\hat{h}^{+}\over \tau}} + {{2\hat{h}^{+}\over \hat{h}^{+} +
\hat{h}^{-}}} \left( {{\hat{h}^{-}\over \hat{h}^{+}}} {\ln } {\hat u_{+} \over
\hat u}  - {{\hat{h}^{+}\over \hat{h}^{-}}} {\ln } {\hat u_{-} \over \hat u}
\right) .
$$
\par
Approximating the system (\ref{sys54}) by these invariants as 
it was done for the semilinear heat equation, we obtain a system
of difference evolution equations. As an example, we present here
an invariant difference model that has explicit equations for the solution
and the trajectory of $x$:
\begin{equation} \label{sh54}
\left\{
\begin{array}{l}
{\Delta x}\ \  =\ \  {\displaystyle {{2\tau\over h^{+} + h^{-}}} }\left( -
{\displaystyle {{h^{-}\over h^{+}}} } {\ln } {\displaystyle {{u_{+}\over u}} }
+ {\displaystyle {{h^{+}\over h^{-}}} } {\ln } {\displaystyle {{u_{-}\over u}}
} \right) ;\\
\\
{\left( {\displaystyle {{u\over \hat{u}}} } \right) }^{2} {\exp } \left( - {
\displaystyle {{1\over 2}} {{ (\Delta x)^{2}\over \tau}} } \right)\ \  =\ \   1
-  {\displaystyle {{4\tau\over h^{+} + h^{-}}} \left( {{1\over h^{+}}} {\ln }
{{u_{+}\over
u}} + {
{1\over h^{-}}} {\ln } {{u_{-}\over u}} \right) }.\\
\end{array}
\right.
\end{equation}

\bigskip

\noindent {\bf 5. Example of an exact solution.}

\medskip
\noindent Let us find the solution of the difference model (\ref{sh54}) for the
heat equation that is invariant with respect to the operator
\begin{equation} \label{oper54}
2 \alpha X_{2} + X_{3}, \ \ \alpha = const.
\end{equation}

The operator (\ref{oper54}) has three invariants: $t$ and the expressions  
$u \exp ( {\textstyle {x^{2}\over 4(t + \alpha )} } )$ and
$({\textstyle {\Delta x\over \tau} - {x\over t + \alpha } } )$. So, we will
seek the invariant solution in the form:
$$
\left\{
\begin{array}{l}
u(x,t)\ \  =\ \  f(t) \exp \left( - {\displaystyle {x^{2}\over 4(t + \alpha )}
} \right)  ;\\
\\
{\displaystyle {\Delta x\over \tau} }\ \  =\ \  g(t) + {\displaystyle {x\over t
+ \alpha } }.\\
\end{array}
\right.
$$
\noindent Substituting this form of the solution to the system (\ref{sh54}),
we obtain ordinary difference equations for $f(t)$ and $g(t)$:
\begin{equation}   \label{eq55}
\left\{
\begin{array}{l}
f(t+\tau)\ \  =\ \  \left( {\displaystyle {t + \alpha \over t + \tau +
\alpha } } \right) ^{1/2} f(t) ;\\
\\
g(t)\ \  =\ \  0 .\\
\end{array}
\right.
\end{equation}
\par
\noindent Solving this system, we find the solution of the
difference equation
$$            
u(x,t) = f(0) \left( {\alpha \over t + \alpha } \right) ^{1/2} \exp \left( -
{x^{2} \over 4(t + \alpha )} \right),
$$            
\noindent and the solution for evolution of a grid:
$$
x = x^{0} \left( {t + \alpha \over \alpha } \right) .
$$
\noindent The obtained solution is the solution of the Cauchy problem
with the invariant initial condition:
$$
u(x,0) = f(0) \exp \left(  - {x^{2}\over 4\alpha } \right) .
$$
\noindent If $\alpha = 0$, the fundamental solution of the heat equation
\begin{equation}  \label{eq56}
u(x,t) = C \left( {1\over t} \right) ^{1/2} \exp \left( - {x^{2} \over 4t }
\right) .
\end{equation}
\noindent is a solution of the difference model. This solution holds on the
grid:
$$
\Delta x = {\tau \over t} x .
$$
\par
In all cases listed above the difference mesh is arbitrary at the initial
point, $t = 0$. In this case it will not be uniform on other time layers.
If the grid is uniform in $x$--direction at $t = 0$
($h_{+} = h_{-} = h$), the steps of the grid in the $x$--direction will be
equal
each other on every time layer, but differ from steps on the previous
time layer.

\par
It is necessary to mention that the obtained difference invariant solution
is the solution of the corresponding differential equation that is
invariant with respect to the operator (\ref{oper54}). As in the 
differential case the above reduction
procedure can be applied for every subalgebra of the algebra (\ref{op54}), 
and then
one obtain different moving meshes which are self-adaptive to every solution.

\par
Thus, the above difference models inherit both the groups of the 
differential equations and the potential to be integrated on a subgroup.

\bigskip

\noindent {\bf 6. Numerical calculations }

\medskip
\noindent
Here we do not discuss the questions of stability and convergence of the
developed
schemes. These are hard questions for nonlinear schemes but one of the
ways to check them is by computing the numerical solutions to the 
exact solutions of the original
differential equations. Below we present the results of numerical calculation
of the invariant solution (\ref{eq56}) by means of the invariant model
(\ref{sh54}).
It is necessary to note that the calculations were not done for the 
equations (\ref{eq55}) reduced on the subgroup, but for the nonstationary
equations (\ref{sh54}). Initial data correspond to the solution
(\ref{eq56}) with $t = 10$. In Fig. 5 we present the evolution of $u$
from invariant initial data by the invariant scheme (\ref{sh54}).


\begin{figure}[h!]    \label{Figinsol}
\centerline{\psfig{figure=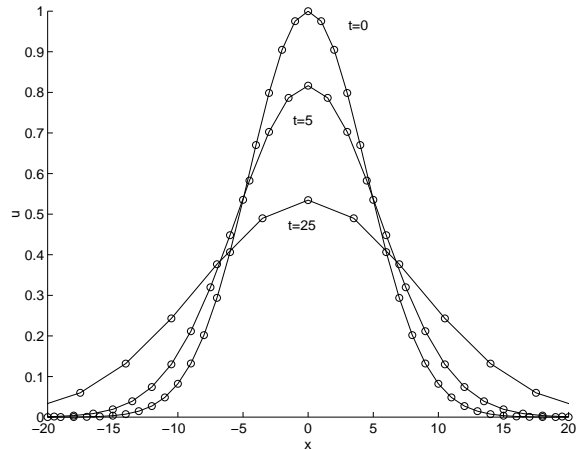,angle=0,width=9cm}}
\caption{Solution of the invariant model.}
\end{figure}

Let us note that the difference model (\ref{sh54}) gives us 
practically the exact solution for equation (\ref{eq54}). 
There is only round-off error of computations.  
In the Fig. 6 the evolution of the grid in the plane $(t,x)$
for the calculation of the solution (\ref{eq56}) is shown.


\begin{figure}[h!]   \label{Figgird}
\centerline{\psfig{figure=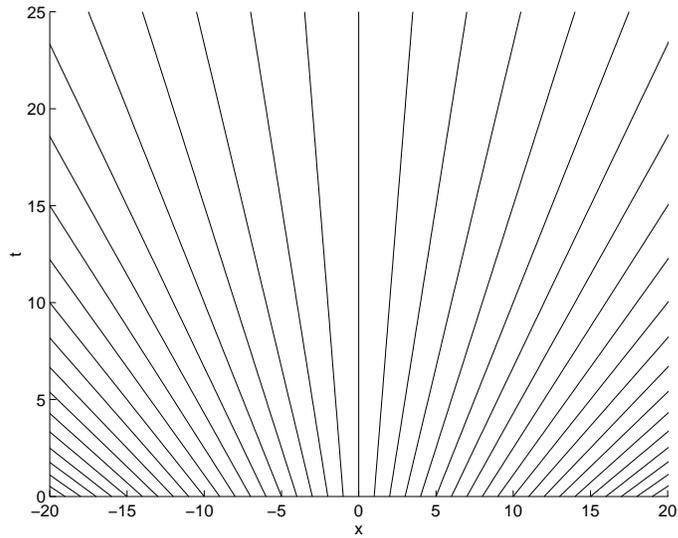,angle=0,width=9cm}}
\caption{Evolution of the mesh.}
\end{figure}



\begin{figure}[h!]
\centerline{\psfig{figure=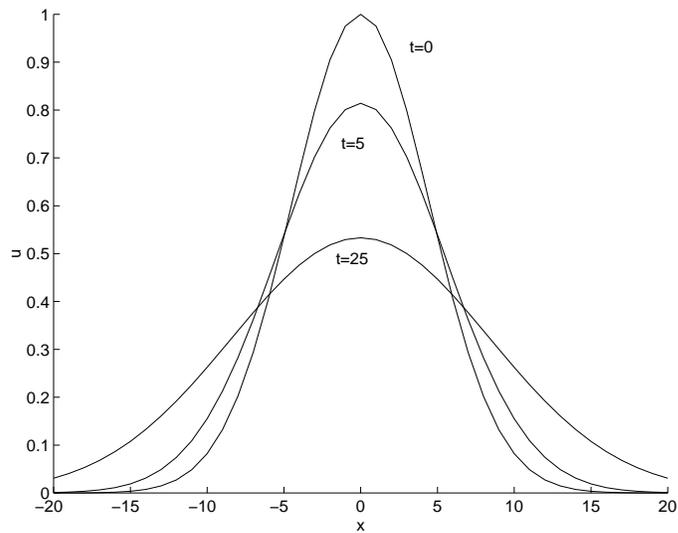,angle=0,width=9cm}}
\caption{Solution of the noninvariant scheme.}
\end{figure}

The same calculations for the difference equation (\ref{eq55a}) on the
orthogonal grid (as that in Fig. 1) gives a similar picture (see Fig. 7). 
In this case the solution does not coincide with exact
solution of
the equation (\ref{eq54}). The difference between exact solution and
numerical solution is shown in Fig. 8.


\begin{figure}[h!]
\centerline{\psfig{figure=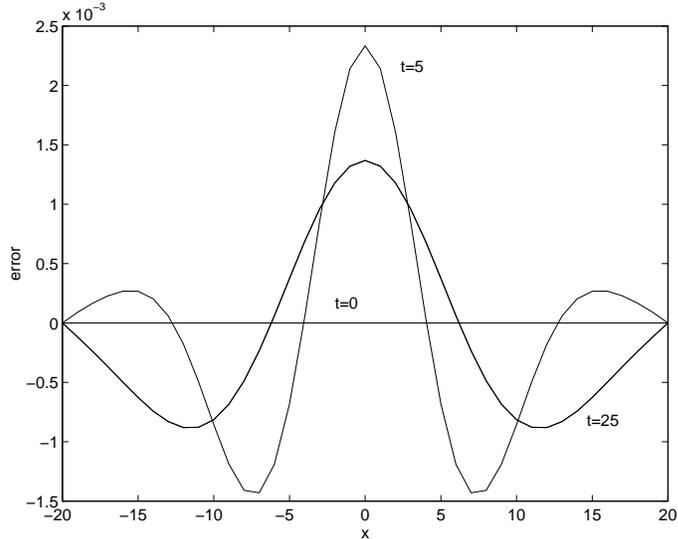,angle=0,width=9cm}}
\caption{The error of the noninvariant scheme.}
\end{figure}

We should note that for the calculation on the model (\ref{sh54}) we
defined the boundary values $u$ on the moving ends of
the space interval. For the difference
equation (\ref{eq55a}) we defined the boundary values of $u$ on the ends of the
fixed orthogonal grid in accordance with the same solution.
The comparison of the two different models shows that for the
invariant model even
on the decreasing number of the points of the grid on the initial
space interval we have greater accuracy than for the noninvariant
one.

\bigskip
\medskip

\noindent {\bf Remarks}

\medskip
\noindent Following the above technique for the Burgers equation for
the potential
$$
w_{t}  + {1 \over 2} w_{x}^{2} = w_{xx},
$$
we obtain the finite-difference model for this equation
on a moving mesh
$$
\left\{
\begin{array}{l}
\Delta x\ \  =\ \  \tau {\displaystyle { h^{-} \dh w { } _{x}  + h^{+}
\dh w { } _{\bar x} \over h^{+} + h^{-} } };\\
\\
\exp \left( \hat{w} - w - { \displaystyle {\Delta x^{2} \over 2\tau} } \right)\
\  =\ \  1 + \tau \dh w { } _{x \bar x} ;\\
\end{array}
\right.
$$
where $\dh w { } _{x \bar{x}}  ={2 \over h^{-} + h^{+} } ( \dh w { }_{x}
- \dh  w { } _{\bar x} )$.
\par
It is
interesting to note that this model is connected with model (\ref{sh54})
by the same Hopf transformation
$$
w = - 2 \ln u 
$$
as their continuous counterparts.

It is important to notice that in all cases the moving in $(x,t)$-
plane meshes can be stopped by the new coordinates of Lagrange's type
with one additional dependent variable (for involving those coordinates see,
for example, [14]).

\bigskip
\noindent { \it  Acknowledgments.} The authors' research was partly
supported
by The Norwegian Research Council under contract no. 111038/410, through 
the SYNODE project, and  Russian Fund for Base Research.

\end{document}